\documentclass[sn-apa]{sn-jnl}

\usepackage[utf8]{inputenc}
\usepackage{bbm}
\usepackage{amsmath}
\usepackage{csquotes}
\usepackage{mathtools}
\usepackage{multirow}
\usepackage{makecell}
\usepackage{amsthm}
\usepackage{adjustbox}
\usepackage{amssymb}
\usepackage{pdflscape}
\usepackage{amsfonts}
\usepackage{anyfontsize}
\usepackage[english]{babel}
\usepackage{xfrac}
\usepackage[hang]{footmisc}
\usepackage{lipsum}
\usepackage{enumitem}
\usepackage{hyperref}
\hypersetup{colorlinks=true, pdfauthor=author, allcolors=blue}
\usepackage{multirow}
\usepackage{graphicx,xcolor,textpos}
\usepackage{emptypage}
\usepackage{appendix}
\usepackage{booktabs}
\usepackage[toc,acronym,nonumberlist]{glossaries}
\usepackage{diagbox}
\usepackage{subcaption}
\usepackage{float}
\usepackage{listings}
\usepackage{lipsum}

\DeclareMathOperator{\EX}{\mathbb{E}}

\newcommand{\indep}{\perp \!\!\! \perp}

\DeclareMathOperator*{\argmax}{arg\,max}

\theoremstyle{thmstyleone}%
\newtheorem{theorem}{Theorem}
\theoremstyle{thmstyletwo}%

\theoremstyle{thmstylethree}%

\begin{document}

\title[Flexible CF approach under dependent censoring]{A flexible control function approach for survival data subject to  different types of censoring}

\author*[1]{\fnm{Ilias} \sur{Willems}}\email{ilias.willems@kuleuven.be}
\equalcont{These authors contributed equally to this work.}

\author[2]{\fnm{Sara} \sur{Rutten}}\email{sara.rutten@uhasselt.be}
\equalcont{These authors contributed equally to this work.}

\author[1]{\fnm{Gilles} \sur{Crommen}}\email{gilles.crommen@kuleuven.be}

\author[1]{\fnm{Ingrid} \sur{Van Keilegom}}\email{ingrid.vankeilegom@kuleuven.be}

\affil*[1]{\orgdiv{ORSTAT}, \orgname{KU Leuven}}
\affil*[2]{\orgdiv{I-BioStat}, \orgname{Hasselt University}}

\abstract{This paper addresses the problem of identifying and estimating the causal effect of a treatment in the presence of unmeasured confounding and various types of right-censoring. Examples of these censoring mechanisms are administrative censoring, competing risks and dependent censoring (e.g. loss to follow-up). Different parametric transformations are applied to each event time, resulting in a regression model with a more additive structure and error terms that are approximately normal and homoscedastic. The transformed event times are modeled using a joint regression framework, assuming multivariate Gaussian error terms with an unspecified covariance matrix. A control function approach is used to deal with unmeasured confounding. The model is shown to be identifiable and a two-step estimation procedure is proposed. This estimator is proven to yield consistent and asymptotically normal estimates. Furthermore, a goodness-of-fit test for the model's validity is developed. Simulations are conducted to examine the finite-sample performance of the proposed estimator under various scenarios. Finally, the methodology is applied to investigate the causal effect of job training programs on unemployment duration using data from the National Job Training Partnership Act (JTPA) study.}


\maketitle


\section{Introduction}

In this paper, we focus on estimating the causal effect of a treatment variable on a right-censored event time. As such, this research is situated at the intersection of causal inference and survival analysis. This is a particularly practicable context since the problems of censoring and unmeasured confounding occur jointly in many fields of applied research, such as medicine or economics. For example, one can be interested in the effect of job training services on unemployment duration. Each participant in the study is assigned to either the treatment or the control group $(W)$, after which each individual can decide whether or not they participate in the job training services $(Z)$. The event time of interest is time until employment $(T)$.

The data resulting from such studies are frequently subject to right-censoring, which can arise from various underlying causes. We will classify these right-censoring events into three categories. The first category consists of censoring events that are stochastically independent of the event of interest. A typical example is administrative censoring, where subjects are censored because they have not found a job before the end of the study. The second type consists of censoring mechanisms that may depend on the event of interest, such as loss to follow-up and drop-out. As an example, a lack of motivation may link early censoring times to long unemployment durations. Assuming that these possibly dependent censoring events are stochastically independent of the event of interest can lead to biased results \citep{emura2018analysis}. Finally, the third category involves competing risks, which prevent the occurrence of the event of interest. These competing events are mutually exclusive and often dependent on each other. For example, leaving the labor force to return to school because of poor employment prospects.

The analysis is complicated further when the treatment is subject to unmeasured confounding. This implies that the causal effect of a treatment variable $Z$ on the event time $T$ cannot be identified from the conditional distribution of $T$ on $(X^{\top},Z)$, where $X$ represents a vector of exogenous covariates. In the example of unemployment duration, it could be that the study participants who \emph{choose} to follow the job training services are those who benefit the most from these services. Therefore, a naive analysis based on the conditional distribution of $T$ on $(X^{\top},Z)$ may result in biased estimates.

To address these challenges, we propose a joint regression model that accounts for all types of right-censoring mentioned before, while also allowing for the identification of the causal effect of the treatment. Each event time, except for administrative censoring, is modeled by a transformed linear regression model where the error terms follow a multivariate Gaussian distribution, allowing for the estimation of the correlations between event times. This transformation will be chosen such that it leads to a more additive structure with approximately normal and homoscedastic error terms. Note that we allow for different transformations for each event time. A control function approach, which relies on an instrumental variable, is leveraged such that the causal effect of the treatment can be identified. We propose a two-step estimation method and prove the consistency and asymptotic normality of the parameter estimates. Based on a weighted Cramér–Von Mises statistic, a goodness-of-fit test is developed to test the validity of the proposed model. The distribution of the test statistic is estimated by a parametric bootstrap procedure. Simulation results under different scenarios and a data application regarding the effect of job training programs on unemployment duration are also provided. All developed code can be found on GitHub at \url{https://github.com/WillemsIlias/FMCC} and is also implemented in the \texttt{depCensoring} package that can be downloaded from CRAN.

\subsection{Related literature} 

A fundamental result in the survival analysis literature comes from \cite{Tsiatis1975}, who proved that the joint distribution of two or more failure times cannot be identified in a nonparametric way by only observing their minimum. This implies that conditions have to be imposed on the joint and/or marginal distributions of the latent survival and censoring times to identify their joint distribution.

The most popular approach to deal with this issue is based on the use of copulas. \cite{ZhengKlein1995} were the first to introduce this idea and proposed the copula-graphic estimator. It extends the well-known \cite{KaplanMeier1958} estimator to the dependent censoring case. Their method allows for a nonparametric estimator of the marginals of $T$ and $C$ under the assumption of a known copula for their joint distribution. In this context, \emph{known copula} means that both the association parameter and functional form of the copula are known. Moreover, \cite{RivestWells2001} studied the copula-graphic estimator in the special case of a known Archimedean copula. The copula-graphic estimator was further extended to account for multiple competing risks by \cite{carriere1995removing}, who demonstrated that the marginal survival functions remain identifiable under this generalization. However, this method involves solving a system of non-linear differential equations, which can become computationally intensive as the number of competing risks increases. To address this issue, \cite{lo2010copula} introduced a more practical approach based on the assumption of an Archimedean copula. Their method simplifies the multiple competing risks model by aggregating irrelevant risks, reducing it to a two-risk model. Therefore, they can directly apply the copula-graphic estimator for Archimedean copulas. This approach was further generalized to the regression setting by \cite{lo2014regression}, allowing for both multiple competing risks and covariates.

As the assumption of a known copula is usually not realistic in practice, an alternative approach for two competing risks was proposed by \cite{Czado2023}. They showed that, in exchange for fully parametric marginals, the association parameter describing the dependence between them can be identified. This can be seen as surprising since we only observe the minimum of the two competing risks. Most recently, \cite{Deresa2023} showed that the association parameter remains identifiable when one of the latent survival times follows a semiparametric \cite{cox1972} proportional hazards model. In subsequent work, \cite{Deresa2024} improved this result to the case where also $C$ follows a semiparametric Cox model. \cite{Deresa2020_a}, \cite{Deresa2020} and \cite{DeresaVKAntonio2022} proposed similar methods that also allow for covariates, multiple competing risks and left truncation, respectively. In the multiple competing risks setting, \cite{Lo2017a} again proposed to pool the irrelevant risks into one variable and model the dependence with the risk of interest using an Archimedean copula. Without specifying a model for the competing risks, they could identify the sign of covariate effects. In subsequent work, \cite{Lo2023} imposed a parametric Cox proportional hazards model for the risk of interest. Moreover, they left the distribution of the pooled irrelevant risks completely unspecified. The authors showed the identifiability of the model and provided an estimator of the association parameter as well as the other finite dimensional parameters. 

Secondly, this research relates to the literature surrounding causal inference, in which instrumental variable methods are commonly used to address the biases induced by unmeasured confounding. In the context of survival analysis, one is often interested in the effect of a treatment on the time until a certain event of interest occurs (cf. Introduction). The confounded variable under investigation is therefore often the binary variable indicating the treatment received. In this context, \cite{richardson2017nonparametric}, \cite{beyhum2023nonparametric} and \cite{martinussen2020instrumental} proposed methodologies to study the treatment effect on several quantities of interest, such as the subdistribution function, its corresponding quantile function, or the cause-specific hazard function of $T$, respectively. However, they avoided certain identifiability issues by assuming censoring to be independent of the competing risks. Under the same independence assumption, \cite{Zheng2017}, \cite{Kjaersgaard2016} and \cite{Ying2019} used instrumental variables to study the subdistribution function of $T$ in a setting where the confounded variable is not imposed to be binary, and \cite{Wang2023} focused on the cause-specific hazard function. However, note that all of the aforementioned papers only do inference on \emph{crude} quantities, like subdistributions or cause-specific hazards, which are identified based on the observed data (in the absence of confounding). Hence, in particular, none of the discussed approaches allow for inference on the marginal distribution of $T$.

To the best of our knowledge, the only available method in the literature that can study the marginal distribution of $T$ in the presence of unobserved confounding and dependent censoring is by \cite{crommen2024instrumental}, who proposed to combine a control function approach with a parametric model for the event and dependent censoring times. Their method, however, imposes a bivariate normality assumption on the error terms, which can be difficult to justify in practice. Moreover, it lacks a specification test. The current paper will extend this work in many meaningful ways. Firstly, the model is generalized to the setting of multiple latent event times, and to additionally allow for administrative censoring, which requires substantial theoretical effort. Secondly, different power transformations are applied to each event time when including them in the model, resulting in a more additive structure with approximately normal and homoscedastic error terms \citep{box1964analysis}. This generalization significantly decreases the stringency of the imposed distributional assumption but brings further theoretical complications. Lastly, we develop a goodness-of-fit test for the overall model proposed in this paper. A data application illustrates the importance of the added generality of our approach over its predecessor.

\subsection{Outline} 

In Section \ref{section:model}, the model is specified and some useful distributions and definitions are given. In Section \ref{section:identandest}, it will be argued that the proposed model is identified, the likelihood will be defined and the estimation method explained. Moreover, we prove the consistency and asymptotic normality of the estimator. In Section \ref{section: Goodness of fit test}, a goodness-of-fit test is proposed. Finally, Section \ref{section: Simulation study} contains a simulation study that assesses the finite sample performance of our methodology and Section \ref{section:dataapp} uses the proposed method to study the causal effect of job training programs on the time until employment.


\section{The model}\label{section:model}

\subsection{Model specification}

Define $T^1, \ldots ,T^K$ to be the logarithm of $K$ latent event times and let $C$ represent the logarithm of a right-censoring time that is independent of $(T^1, \ldots, T^K)$. A common case corresponds to $K = 2$, where $T^1$ represents the event time of interest and $T^2$ represents the dependent censoring time. Define the observed follow-up time $Y=\min(T^1, \ldots, T^K,C)$ and the indicator $\Tilde{\Delta}=(\Delta_1, \ldots ,\Delta_K,\Delta_C)$ with $\Delta_k=\mathbbm{1}(Y=T^k)$ for $k=1,\ldots, K$ and $\Delta_C=\mathbbm{1}(Y=C)$. Note that the following model specification is also valid when there is no administrative censoring $(C\equiv \infty)$. The exogenous covariates are given by $X=(1,\tilde{X}^{\top})^{\top}$, where $\tilde{X}$ is of dimension $m$. We also allow for a possibly endogenous variable $Z$, which is assumed to be univariate. More precisely, we propose the following joint regression model:
\begin{equation}
	\label{eq:Regression Model rewritten}
		\Lambda_{\theta_k}(T^k)=X^{\top}\beta_k+Z\alpha_k+V\lambda_k+\epsilon_k, \quad k=1,\ldots, K,
\end{equation}
where $(\epsilon_1, \ldots, \epsilon_K)$ are unobserved error terms, $\Lambda_{\theta_k}(\cdot)$ is a power transformation known up to $\theta_k$, such that $(\theta_1,\dots\theta_K)$ are all allowed to be different from each other, and $V$ is an unobserved confounder of $Z$. Since $V$ is not observed, we use an instrumental variable $\tilde{W}$ that is sufficiently dependent on $Z$ (conditionally on $X$) to be able to identify the causal effect of interest $\alpha_{k}$. More precisely, we let $V=g_\gamma(W,Z)$, with $W=(X^{\top},\Tilde{W})$, for which the control function $g$ is known up to the parameter $\gamma$. Note that this control function follows from the reduced form, which is specified by the analyst. A simple example can be given by
\begin{equation}
    V  = Z-W^\top \gamma, \label{eq:CF_Z_continuous}
    \end{equation}
    when $Z$ is continuous and linearly related to $W$, that is,  $$Z = W^\top \gamma + \nu \quad \text{with} \quad \EX[\nu \mid W] = 0.$$ Another, more involved, example of a control function is
    \begin{equation}
    V  = Z\EX[\nu   \mid  W^\top{\gamma} > \nu ] + (1-Z)\EX[\nu   \mid   W^\top{\gamma} < \nu ], \label{eq:CF_Z_binary}     
    \end{equation}
when $Z$ is binary and the relation between $Z$ and $W$ is specified as $$Z = \mathbbm{1}(W^\top \gamma > \nu) \quad \text{with} \quad \nu \indep W.$$ 
We refer to \cite{wooldridge2010econometric}, \cite{Navarro2010} and \cite{Tchetgen2015} for further justification and other examples of the control function approach. Moreover, it is assumed that:
\begin{enumerate}[label={(A\arabic*)}]
\item \label{as:A1} 
    $(\epsilon_1, \dots, \epsilon_K)^\top
    \sim \mathcal{N}_{K}\left(
            \textbf{0 }, \Sigma\right),$ where $\Sigma \in \mathbb{R}^{K \times K}$ is a positive definite matrix.
    \item \label{as:A2} $(\epsilon_1 ,\ldots ,\epsilon_K,C) \indep (W^{\top},Z)$.
    \item \label{as:A3} $(T^1 ,\ldots ,T^K)$ and $C$ are conditionally independent, given $(W^{\top},Z)$.
    \item \label{as:A4}The covariance matrix of $(\tilde{X}^{\top},Z,V)$ has full rank and $\text{Var}(\tilde{W}) > 0.$
    \item \label{as:A5} There exists an $M_Y \in \mathbb{R}$, specified in Appendix C.1, such that for all $k \in \{1,\dots,K\}$ the probabilities $P(\Delta_k = 1 \vert Y=y, W=w,Z=z)$ are strictly positive for almost all $(W^{\top},Z)$ and for all $ y \in (-\infty,M_Y)$.
    \item \label{as:A6} Censoring by $C$ is non-informative for $(T^1, \ldots ,T^K)$ given $(W^{\top},Z)$.
    \item \label{as:A7} $\{\Lambda_\theta:\theta \in \Theta\}$ is a family of strictly increasing, continuously differentiable transformations, defined on the whole real line and for which it holds that $\lim_{t \rightarrow a}{\Lambda_\theta(t)} = a$ for all $\theta \in \Theta$ with $a \in \{-\infty,+\infty\}$.
    \item \label{as:AX} For all $\theta^1, \theta^2 \in \Theta$ it holds that $\lim_{y \to -\infty}\Lambda_{\theta^1}(y) / \Lambda_{\theta^2}(y) = \infty \iff \theta^1 < \theta^2$.
    \item \label{as:A8} For all $\theta^j \in \Theta, \mu^j \in \mathbb{R}, \sigma^j > 0$, $j=1,2$, it holds that if $\theta^1 \neq \theta^2$, then the limit $$\lim_{t \rightarrow \pm \infty} {K_{\theta^1,\mu^1,\sigma^1}(t) /K_{\theta^2,\mu^2,\sigma^2}(t) =0 \text{ or }\infty},$$ with 
\begin{equation*}
K_{\theta,\mu,\sigma} (t) = \exp \biggl(-\frac{1}{2} \biggl\{ \frac{\Lambda_\theta(t)-\mu}{\sigma} \biggl \}^2 \biggl) \Lambda_\theta'(t).
\end{equation*}
\end{enumerate}
From Assumption \ref{as:A1}, it is clear that that the transformed survival times $(\Lambda_{\theta_1}(T^1),\dots,\Lambda_{\theta_K}(T^K))$, conditional on $(X^{\top},Z,V)$, are multivariate normally distributed and allowed to be dependent on each other. This is a strong assumption to make, but it will later be shown by Theorem \ref{identificationtheorem} that this allows us to actually identify the variance-covariance matrix $\Sigma$. This can be seen as surprising, as it means that we can identify the relationship between $(\Lambda_{\theta_1}(T^1),\dots,\Lambda_{\theta_K}(T^K))$ while only observing the minimum of $(T^1, \ldots, T^K,C)$. Note that Assumption \ref{as:A1} is placed on the error terms of the transformed survival times. Therefore, $\Lambda_{\theta}$ should be a transformation that aims to improve normality. It is well known that applying a power transformation to the response of a regression model results in a more additive structure with approximately normal and homoscedastic error terms \citep{box1964analysis}. The most commonly used example of such a power transformation is the Box-Cox transformation. However, it is already clear that this transformation does not satisfy Assumption \ref{as:A7}. A transformation that does satisfy this assumption is the Yeo-Johnson transformation, which can be seen as an extension of the Box-Cox transformation to the whole real line \citep{yeo2000new}. This transformation can be expressed as:
\begin{equation}
	\label{eq:Yeo-Johnson}
	\Lambda_\theta(t)=
	\left\{
	\begin{array}{ll}
		\{(t+1)^\theta-1\}/\theta & t\geq 0, \theta \neq 0 \\
		\log(t+1) & t\geq 0, \theta = 0 \\
		-\{(-t+1)^{2-\theta}-1\}/(2-\theta) & t < 0, \theta \neq 2 \\
		-\log(-t+1) & t< 0, \theta = 2 \\
	\end{array}
	\right..
\end{equation}
Note that when $\theta = 1$, this transformation is the identity transformation. When $1<\theta \leq 2$, the transformation is convex, implying a contraction of the lower part of the support and an extension of the upper part, decreasing skewness to the left. When $0 \leq \theta <1$, the transformation is concave and decreases the skewness to the right. Hence, the transformation tries to improve the symmetry of the distribution. As is the case for the Box-Cox transformation, the Yeo-Johnson transformation aims to improve normality, making Assumption \ref{as:A1} more plausible. It is clear that Assumption \ref{as:AX} is satisfied by the Yeo-Johnson transformation by direct calculation of the limits. Section A of the Supplementary Material verifies that this transformation satisfies Assumption \ref{as:A8}.
 
\subsection{Useful definitions and distributions}

For the remainder of this paper, it will be useful to have some definitions and distributions at hand so we can more easily state and prove our results. We observe the random vector $S = (Y,\Tilde{\Delta},\tilde{X},\tilde{W},Z)$, which takes values in the space $\mathcal{G} = \mathbb{R} \times \{0,1\}^{K+1} \times \mathbb{R}^{m} \times \mathbb{R} \times \mathbb{R}$. Furthermore, denote the parameter space of $\gamma$ by $\Gamma \subset \mathbb{R}^{m+2}$ and collect all other parameters in the vector $\eta$. The parameter space of $\eta$ is given by $\mathcal{H} \subset \{ \eta: (\beta_1^{\top}, \ldots, \beta_K^{\top},\alpha_1, \ldots ,\alpha_K, \lambda_1 ,\ldots ,\lambda_K) \in \mathbb{R}^{K(m+3)},(\sigma_1 ,\ldots ,\sigma_K) \in \mathbb{R}^K_{> 0},(\rho_{12}, \ldots, \rho_{(K-1)K})\in (-1,1)^{K(K-1)/2},(\theta_1, \ldots, \theta_K) \in \Theta^K \}$, with $\Theta = [0,2]$, $(\sigma_1, \dots, \sigma_K)$ represent the square root of the diagonal elements in the covariance matrix $\Sigma$ and $(\rho_{12}, \ldots, \rho_{(K-1)K})$ denote the off-diagonal elements of the corresponding correlation matrix $\Omega$. To ease notation, define: $$b_{k}=\Lambda_{\theta_k}(y)-x^{\top}\beta_k-z\alpha_k-g_\gamma(w,z)\lambda_k,$$ and $$\tau_k=x^{\top}\beta_k+z\alpha_k+g_\gamma(w,z)\lambda_k.$$ Furthermore, let $\bar{\Phi}(y_1,y_2,\dots,y_{K-1};\Sigma) = \Phi(-y_1,-y_2,\dots,-y_{K-1};\Sigma) $ with $\Phi(\cdot;\Sigma)$ a $(K-1)$-variate normal distribution with mean equal to zero, and covariance matrix $\Sigma$. In Section B of the supplementary Material, we derive that
\begin{align}
\begin{split}
   & f_{Y,\Tilde{\Delta}\vert W,Z}(y,1,0,\ldots ,0\vert w,z,\gamma;\eta) \\[5pt] 
   & \quad = \frac{1}{\sigma_1}\bar{\Phi}\biggl(\frac{\Lambda_{\theta_2}(y)-m_{2.1}}{\xi_{2.1}},\ldots, \frac{\Lambda_{\theta_K}(y)-m_{K.1}}{\xi_{K.1}};\Omega_1\biggl)\phi\biggl(\frac{b_{1}}{\sigma_1}\biggl) \\
   & \qquad \times \Lambda_{\theta_1}'(y)P(C>y), \notag
\end{split}
\end{align}
with
$\Omega_1$ a correlation matrix with $(\rho_{23.1}, \ldots, \rho_{(K-1)K.1})$ denoting the off-diagonal elements, defined as 
\begin{equation*}
    \rho_{lk.1}=\frac{\rho_{lk}-\rho_{1l}\rho_{1k}}{\biggl\{ (1-\rho_{1l}^2)(1-\rho_{1k}^2)\biggl\}^{1/2}},
\end{equation*} for $2\leq l<k\leq r$ and where for $k=2,\ldots ,K$:
\begin{equation*}
m_{k.1}=\tau_k+\rho_{1k}\frac{\sigma_k}{\sigma_1}b_1 \quad \text{and} \quad    \xi_{k.1}=\sigma_j(1-\rho_{1k}^2)^{1/2}.
\end{equation*}
 The formulas for the sub-densities of $(Y,\Tilde{\Delta})$ given $(W,Z)$ when $\Delta_k=1$ for $k=2,\ldots, K$ can be derived similarly. Moreover, when $C$ is observed, meaning $\Delta_C=1$, it can be derived that:
\begin{align*}
   & f_{Y,\Tilde{\Delta}\vert W,Z}(y,0,\ldots ,0,1\vert w,z,\gamma;\eta) = \bar{\Phi}\biggl(b_{1}, \ldots,b_{K}; \Sigma\biggl) f_C(y).
\end{align*}

\noindent The log-likelihood of our model can now be defined as:
\begin{equation*}
    l: \mathcal{G} \times \Gamma \times \mathcal{H} \rightarrow \mathbb{R}:(s,\gamma,\eta) \rightarrow l(s,\gamma,\eta)=\log f_{Y,\Tilde{\Delta} \vert W,Z}(y,\Tilde{\delta} \vert w,z,\gamma;\eta),
\end{equation*}
such that the expected log-likelihood is equal to:
\begin{equation}
    \label{eq: expected loglikelihood}
    L(\gamma,\eta)=\EX[l(S,\gamma,\eta)] = \int_{\mathcal{G}}l(s,\gamma,\eta)dG(s),
\end{equation}
where $G$ is the distribution function of $S$.

In the framework of competing risks, the cumulative incidence function (CIF) is often calculated. In the presence of $K$ competing events, the function is defined by $I_k(y\vert W,Z)=P(Y\leq y, \Delta_k=1 \vert W,Z)$ for $k=1,\ldots,K$. Hence, it can be interpreted as the probability of failing from cause $k$ before time $y$ in the presence of the other competing events. For the first cause of failure, this becomes:
 \begin{align}
 \begin{split}
   & I_1(y \vert w,z) \\
    & \quad =\int_{-\infty}^{y}{}\frac{1}{\sigma_1}\bar{\Phi}\biggl(\frac{\Lambda_{\theta_2}(u)-m_{2.1}}{\xi_{2.1}},\ldots, \frac{\Lambda_{\theta_K}(u)-m_{K.1}}{\xi_{K.1}};\Omega_1\biggl) \notag\\
    & \qquad \times \phi\biggl(\frac{\Lambda_{\theta_1}(u)-\tau_1}{\sigma_1}\biggl)\Lambda_{\theta_1}'(u)P(C>u)du.
\end{split}
\end{align}
The cumulative incidence functions for the other causes can be derived similarly.

In the framework of dependent censoring, where $K=2$, we are interested in the marginal distribution of one of the latent event times, say $T^1$.  From Assumptions \ref{as:A1} and \ref{as:A2}, it follows that 
\begin{equation*}
    \label{eq:distribution T}
     F_{T^1 \vert W,Z}(t \vert w,z,\gamma; \eta) = P(T^1\leq t \vert w,z) =\Phi \biggl(\frac{b_1}{\sigma_1} \biggl).
\end{equation*}
Note how, compared to the CIF, the marginal distribution is better suited to study $T^1$ when $T^2$ is a non-terminal, nuisance event. 

\section{Model identification and estimation}\label{section:identandest}

\subsection{Model identification}
Even though we only observe the minimum of the event times and $C$ through the follow-up time $Y$ and censoring indicators $\Tilde{\Delta}$, it can be shown that model \eqref{eq:Regression Model rewritten} is identified. With identified, it is meant that two different sets of the parameters $(\gamma, \eta)$ imply two different joint distributions of $S$. Suppose that $(\gamma^*,\eta^*)$ are the true values of the parameters. We will assume that
\begin{enumerate}[label=(A\arabic*),resume]
\item \label{as:A9} $\gamma^*$ is identified.
\end{enumerate}
It is clear that Assumption \ref{as:A9} holds by Assumption \ref{as:A4} when $Z$ is a continuous random variable for which \eqref{eq:CF_Z_continuous} holds. Moreover, Assumption \ref{as:A9} is implied by Assumption \ref{as:A4} and a known distributional assumption on $\nu$ when $Z$ is a binary random variable for which \eqref{eq:CF_Z_binary} holds \citep{manski1988identification}. We now state our identifiability theorem, for which the proof can be found in Section C.1 of the Supplementary Material:

\begin{theorem}\label{identificationtheorem}
Under Assumptions \ref{as:A1} to \ref{as:A9}, suppose that $(T^1_{j}, \ldots, T^K_{j},C_j)$, with $j=1,2$, satisfy model \eqref{eq:Regression Model rewritten} with parameter vectors $(\gamma,\eta_1)$ and $(\gamma,\eta_2)$, respectively. Denote, $Y_j=\min(T^1_{j}, \ldots,T^K_{j},C_j)$ and
$\Tilde{\Delta}_j = \bigl(\mathbbm{1}(Y_j=T^1_{j}), \ldots,$
$ \mathbbm{1}(Y_j=T^K_{j}), \mathbbm{1}(Y_j=C_j)\bigl)$. Then if $$f_{Y_1,\tilde{\Delta}_1\vert W,Z}(\cdot,l_1, \ldots, l_K,l_C\vert w,z,\gamma;\eta_1)=f_{Y_2,\tilde{\Delta}_2\vert W,Z}(\cdot,l_1, \ldots, l_K,l_C\vert w,z,\gamma;\eta_2)$$ for $l_1, \ldots, l_K,l_C  \in \{0,1\}$ and for almost every $(w,z)$, it follows that:
\begin{equation}
\eta_1 = \eta_2.
\end{equation}
\end{theorem}

It can be noted that \eqref{eq:Regression Model rewritten} is a parametric copula model. Indeed, we impose that parametric transformations of $T^1, \ldots, T^K$ follow a certain normal distribution and model the dependence structure with a Gaussian copula. This makes model \eqref{eq:Regression Model rewritten} similar to the one of \cite{Czado2023}, who study a broader family of models by allowing for different choices of the marginals and copula but do not include covariates in their analysis (and a fortiori do not need to treat the issue of endogeneity). As is also the case here, the identification of such models proves to be the main difficulty. Moreover, the additional flexibility introduced by allowing different transformations for $T^1, \ldots, T^K$ further complicates the identifiability proof, as it is possible that $T^1 < T^2$ but $\Lambda_{\theta_1}(T^1) > \Lambda_{\theta_2}(T^2)$.

\subsection{Estimation of the model parameters}
To estimate the model parameters, we use a two-step procedure which we refer to as the \emph{two-step estimator}. In the first step, $\gamma$ is estimated. Its estimate, $\hat{\gamma}$, can then be used in the second step to estimate $\eta$. Note that using $\hat{\gamma}$ instead of the true (unknown) value will increase the variance of the estimator for $\eta$. For the first step in the estimation procedure, we make the following assumptions:
\begin{enumerate}[label=(A\arabic*),resume] 
\item \label{as:A9.5} The data consist of i.i.d. observations $\{Y_i,\Tilde{\Delta}_i,W_i,Z_i\}_{i = 1, \dots, n}$.
\item \label{as:A10} There exist a known function $m:(w,z,\gamma) \in \mathbb{R}^{m+2} \times \mathbb{R} \times \Gamma \rightarrow m(w,z,\gamma)$ that is twice continuously differentiable with respect to $\gamma$ such that 
\begin{equation*}
\hat{\gamma} \in \argmax_{\gamma \in \Gamma}{n^{-1}\sum_{i=1}^{n}{m(W_i,Z_i,\gamma)}}
\end{equation*}
is a consistent estimator for the true parameter $\gamma^*$.
\end{enumerate}
In the case of control function $\eqref{eq:CF_Z_continuous}$, it can easily be shown that Assumption \ref{as:A10} follows directly from Assumption \ref{as:A4} by using ordinary least squares to estimate $\gamma$. When $Z$ is a binary random variable for which control function \eqref{eq:CF_Z_binary} is specified, maximum likelihood estimation can be used to consistently estimate $\gamma$ under Assumption \ref{as:A4} and a known distributional assumption on $\nu$ \citep{AldrichNelson1986}.

With the estimate $\hat{\gamma}$ from the first step, we can continue with estimating $\eta$ in the second step. Based on the definition of the expected log-likelihood in \eqref{eq: expected loglikelihood}, the log-likelihood function that will be maximized as a function of $\eta$ over the parameter space $\mathcal{H}$ is given by:
\begin{align}
    \label{eq: likelihood competing risks}
    &\hat{L}(\hat{\gamma},\eta)=\prod_{i=1}^{n}\biggl\{ \prod_{j=1}^{K}\bar{f}_{Y,\Tilde{\Delta}\vert W,Z}(Y_i,\Tilde{\Delta}_i \vert W_i,Z_i,\hat{\gamma};\eta)^{\Delta_{j,i}}\biggl \} \bar{f}_{Y,\Tilde{\Delta}\vert W,Z}(Y_i,\Tilde{\Delta}_i \vert W_i,Z_i,\hat{\gamma};\eta)^{\delta_{C,i}},
\end{align}
with $\bar{f}_{Y,\Tilde{\Delta}\vert W,Z}(Y_i,\Tilde{\Delta}_i \vert W_i,Z_i,\hat{\gamma};\eta)$ being $f_{Y,\Tilde{\Delta}\vert W,Z}(Y_i,\Tilde{\Delta}_i \vert w_i,z_i,\hat{\gamma};\eta)$ without the density and distribution function of $C$ since the administrative censoring is assumed to be non-informative (Assumption \ref{as:A6}). The estimator $\hat{\eta}$ is then equal to:
\begin{equation}
\label{eq: MLE competing risk}
\begin{split}
    \hat{\eta}&=(\hat{\beta}_1^{\top}, \ldots, \hat{\beta}_K^{\top},\hat{\alpha}_1, \ldots, \hat{\alpha}_K, \hat{\lambda}_1 ,\ldots ,\hat{\lambda}_K,\hat{\sigma}_1, \ldots, \hat{\sigma}_K, \hat{\rho}_{12}, \ldots ,\hat{\rho}_{(K-1)K},\hat{\theta}_1, \ldots, \hat{\theta}_K)^{\top} \\ & \in \argmax_{\eta \in \mathcal{H}} \hat{L}(\hat{\gamma},\eta).
\end{split}
\end{equation}

\subsection{Consistency and asymptotic normality}
\label{section: Consistency and asymptotic normality}
First, some notation that will be useful in this section is defined:
\begin{align*}
&h_l(S,\gamma,\eta)=\nabla_{\eta}l(S,\gamma,\eta), & &H_{\eta}=\EX[\nabla_{\eta}h_l(S,\gamma^*,\eta^*)],\\
&h_m(W,Z,\gamma)=\nabla_{\gamma}m(W,Z,\gamma), & &H_\gamma = \EX[\nabla_\gamma h_l(S,\gamma^*,\eta^*)], \\
&M=\EX[\nabla_\gamma h_m(W,Z,\gamma^*)], & &\Psi = -M^{-1}h_m(W,Z,\gamma^*), \\
&\tilde{h}(S,\gamma,\eta)=(h_m(W,Z,\gamma)^{\top},h_l(S,\gamma,\eta)^{\top})^{\top}, & &H= \EX[\nabla_{\gamma,\eta}\tilde{h}(S,\gamma^*,\eta^*)].
\end{align*}
Assumption \ref{as:A10} already states that $\hat{\gamma}$ is a consistent estimator for $\gamma^*$. Under some additional assumptions, it can also be proved that $\hat{\eta}$ is a consistent and asymptotically normal estimator of $\eta^*$. These assumptions are:
\begin{enumerate}[resume,label=(A\arabic*)]
\item \label{as:A11} The parameter space $\mathcal{H}$ is compact and $\eta^*$ is an interior point.
\item \label{as:A12} A function $\mathcal{D}(s)$, integrable with respect to $G$, and a compact neighbourhood $\mathcal{N}_\gamma \subset \Gamma$ of $\gamma^*$ exist such that $\vert l(s,\gamma,\eta)\vert \leq \mathcal{D}(s)$ for all $\gamma \in \mathcal{N}_\gamma$ and $\eta \in \mathcal{H}$.
\item \label{as:A13} $\EX\left[\Vert \tilde{h}(S,\gamma^*,\eta^*)\Vert^2\right] < \infty$ and $\EX\left[\sup_{(\gamma,\eta)\in \mathcal{N}_{\gamma,\eta}}\Vert \nabla_{\gamma,\eta} \tilde{h}(S,\gamma,\eta)\Vert \right] < \infty$ where $\mathcal{N}_{\gamma,\eta}$ is a neighbourhood of $(\gamma^*,\eta^*)$ in $\Gamma \times \mathcal{H}$.
\item \label{as:A14} The matrix $H^{\top}H$ is nonsingular.
\end{enumerate}
As usual, $\Vert \cdot \Vert$ denotes the Euclidean norm. Assumptions \ref{as:A11}, \ref{as:A13} and \ref{as:A14} are commonly made assumptions in maximum likelihood theory. Assumption \ref{as:A12} imposes the existence of a dominating function for the log-likelihood over a suitable subset of its domain. The following two theorems can now be formulated:
\begin{theorem}
\label{consistency theorem}
Assuming \ref{as:A1} until \ref{as:A12} and $\hat{\eta}$ as defined in \eqref{eq: MLE competing risk}, it holds that
\begin{equation*}
\hat{\eta} \xrightarrow{p} \eta^*.
\end{equation*}
\end{theorem}
\begin{theorem}
\label{normality theorem}
Assuming \ref{as:A1} until \ref{as:A14} and $\hat{\eta}$ as defined in \eqref{eq: MLE competing risk}, it holds that
\begin{equation*}
\sqrt{n} (\hat{\eta} - \eta^*) \xrightarrow{d} \mathcal{N}(0,\Sigma_\eta),
\end{equation*}
with
\begin{equation*}
\Sigma_\eta=H_\eta^{-1}\EX[\{h_l(S,\gamma^*,\eta^*)+H_\gamma \Psi\} \{ h_l(S,\gamma^*,\eta^*)+H_\gamma \Psi\}^{\top}](H_\eta^{-1})^{\top}.
\end{equation*}
\end{theorem}
The proofs of these two theorems can be found in Sections C.2 and C.3 of the Supplementary Material respectively. 

\subsection{Asymptotic variance} \label{section: Asymptotic variance}
Following Theorem \ref{normality theorem},
a consistent estimator $\hat{\Sigma}_\eta$ for the covariance matrix can be defined as:
\begin{equation*}
\label{eq: Asymptotic variance}
\hat{\Sigma}_\eta=\hat{H}_\eta^{-1}[n^{-1}\sum_{i=1}^{n}{\{ h_l(S_i,\hat{\gamma},\hat{\eta})+\hat{H}_\gamma\hat{\Psi}_i\} \{h_l(S_i,\hat{\gamma},\hat{\eta})+\hat{H}_\gamma\hat{\Psi}_i \}^{\top}}](\hat{H}_\eta^{-1})^{\top},
\end{equation*}
with
\begin{align*}
 &\hat{H}_{\eta}=n^{-1}\sum_{i=1}^{n}{\nabla_{\eta}h_l(S_i,\hat{\gamma},\hat{\eta})}, & &\hat{H}_\gamma = n^{-1} \sum_{i=1}^{n}{\nabla_\gamma h_l(S_i,\hat{\gamma},\hat{\eta})}, \\
&\hat{M}=n^{-1}\sum_{i=1}^{n}{\nabla_\gamma h_m(W_i,Z_i,\hat{\gamma})}, & &\hat{\Psi}_i = -\hat{M}^{-1}h_m(W_i,Z_i,\hat{\gamma}).
\end{align*}
Using this asymptotic covariance matrix estimate, asymptotic confidence intervals for $\eta$ can easily be constructed. To avoid negative values in the confidence interval of $\sigma_k$, one could construct them on the logarithmic scale using the Delta method and then transform them back to the original scale. Similarly, a Fisher's z-transformation can be used to construct confidence intervals for $\rho_{kj} \in (-1,1)$.

\section{Goodness-of-fit test}
\label{section: Goodness of fit test}
In the following section, we discuss a goodness-of-fit test for the special case of model \eqref{eq:Regression Model rewritten} with only two competing risks $(K=2)$ for ease of explanation. However, the test can be extended to the case of $K>2$ latent event times. The observed times $Y$ correspond to events $(T^1, T^2)$ or administrative censoring $(C)$. However, since we only model $T^1$ and $T^2$, we will not construct a goodness-of-fit test for all observed times $Y = \min(T^1,T^2, C)$, but rather define $\tilde{T} = \min(T^1,T^2)$ and construct a goodness-of-fit test based on $\tilde{T}$, which will be right-censored by $C$. The test will be based on a weighted Cramér--Von Mises type statistic, for which we will estimate the distribution with a parametric bootstrap procedure (see \cite{ET1998} for more details). The test will be such that rejection of the null hypothesis indicates a bad fit, but failing to reject the null hypothesis does not necessarily indicate a good fit. This is because, like most statistical tests, there is the possibility of making a type-I error, but additionally, we are only testing on $\tilde{T}$ and not on $T^1$ and $T^2$ separately. This means that we cannot directly test the goodness-of-fit of the models for $T^1$ and $T^2$ separately, but can only gauge it through testing the goodness-of-fit of the overall model. The reason for this caveat is that directly testing on $T^1$ and $T^2$ is not possible due to their nonparametric joint model being unidentified \citep{Tsiatis1975}. To make this more precise, the null hypothesis of the test is
\begin{equation*}
    H_0: P(\tilde{T} \leq \tilde{t}) = F_{\tilde{T}}(\tilde{t}; \gamma^*, \eta^*) \quad \forall \tilde{t} \in \mathbb{R},
\end{equation*}
where $F_{\tilde{T}}(\tilde{t}; \gamma^*, \eta^*)$ denotes the cumulative distribution function of $\tilde{T}$ that follows from model \eqref{eq:Regression Model rewritten} with $K=2$ using the true set of parameters $(\gamma^*, \eta^*)$. It further follows that
\begin{align*}
    F_{\tilde{T}}(\tilde{t}; \gamma, \eta) & = \iiint \Phi\left(\frac{\Lambda_{\theta_1}(\tilde{t}) - \tau_1}{\sigma_1}\right)f_{\tilde{X}, \tilde{W}, Z}(x,w,z) dxdwdz \\
    & \quad + \iiint \Phi\left(\frac{\Lambda_{\theta_2}(\tilde{t}) - \tau_2}{\sigma_2}\right)f_{\tilde{X}, \tilde{W}, Z}(x,w,z) dxdwdz\\
    & \quad - \iiint \Phi\left(\frac{\Lambda_{\theta_1}(\tilde{t}) - \tau_1}{\sigma_1}, \frac{\Lambda_{\theta_2}(\tilde{t}) - \tau_2}{\sigma_2}; \rho_{12}\right)f_{\tilde{X}, \tilde{W}, Z}(x,w,z) dxdwdz,
\end{align*}
where $f_{\tilde{X}, \tilde{W}, Z}(x,w,z)$ is the joint density of $\tilde{X}$, $\tilde{W}$ and $Z$. We omitted the integration region $\mathbb{R}^m \times \mathbb{R} \times \mathbb{R}$ from the notation of the integral. From this last expression, an estimator of $F_{\tilde{T}}(\tilde{t}; \gamma^*, \eta^*)$ is obtained. Suppose $(\hat{\gamma}, \hat{\eta})$ are the estimated parameters of the model and  define $\hat{\tau}_{1, i} = x_i^{\top}\hat{\beta}_1 + z_i\hat{\alpha}_1 + \hat{v}_i\hat{\lambda}_1$ and $\hat{\tau}_{2, i} = x_i^{\top}\hat{\beta}_2 + z_i\hat{\alpha}_2 + \hat{v}_i\hat{\lambda}_2$. Then we can estimate $F_{\tilde{T}}(\tilde{t}; \gamma^*, \eta^*)$ by
\begin{align*}
    \hat{F}_{\tilde{T}}(\tilde{t}; \hat{\gamma}, \hat{\eta}) & = \frac{1}{n}\sum_{i=1}^n \Phi\left(\frac{\Lambda_{\hat{\theta}_1}(\tilde{t}) - \hat{\tau}_{1, i}}{\hat{\sigma}_1}\right) + \frac{1}{n}\sum_{i=1}^n \Phi\left(\frac{\Lambda_{\hat{\theta}_2}(\tilde{t}) - \hat{\tau}_{2, i}}{\hat{\sigma}_2}\right)\\
    & \quad - \frac{1}{n}\sum_{i=1}^n \Phi\left(\frac{\Lambda_{\hat{\theta}_1}(\tilde{t}) - \hat{\tau}_{1, i}}{\hat{\sigma}_1}, \frac{\Lambda_{\hat{\theta}_2}(\tilde{t}) - \hat{\tau}_{2, i}}{\hat{\sigma}_2}; \hat{\rho}_{12}\right).
\end{align*}
Furthermore, since Assumptions \ref{as:A2} and \ref{as:A3} together imply that $\tilde{T}$ and $C$ are independent of each other, unconditionally on the value of the covariates $(W, Z)$, we can estimate $P(\tilde{T} \leq \tilde{t})$ based on the \cite{KaplanMeier1958} estimator
\begin{equation*}
    \hat{F}_{\tilde{T}, n}(\tilde{t}) = 1 - \prod_{i: \tilde{T}_{(i)} \leq \tilde{t}} \left(1 - \frac{d_{(i)}}{N_{(i)}}\right),
\end{equation*}
where $\tilde{T}_{(i)}$ is the $i$-th ordered observation of $\tilde{T}$, $N_{(i)} = \sum_{j=1}^n I(Y_j \geq \tilde{T}_{(i)})$ and $d_{(i)} = \sum_{j = 1}^n I(Y_j = \tilde{T}_{(i)}, \Delta_{1,j} + \Delta_{2,j} = 1))$. The null hypothesis can then be tested using the following weighted Cramér--Von Mises type statistic:
\begin{equation*}
    T_{CM} = n\int_{\mathbb{R}} \left(\hat{F}_{\tilde{T}}(\tilde{t}; \hat{\gamma}, \hat{\eta}) - \hat{F}_{\tilde{T}, n}(\tilde{t})\right)^2 w(\tilde{t}) d \hat{F}_{\tilde{T}}(\tilde{t}; \hat{\gamma}, \hat{\eta}).
\end{equation*}
The weight function $w(\tilde{t})$ is often taken equal to $1$ everywhere when doing this test based on the parametric bootstrap procedure explained below. In general, large values of $T_{CM}$ will indicate that the model is misspecified. More precisely, we will reject $H_0$ on a $(1 - \kappa)$-confidence level if $T_{CM}$ is larger than the $100 \times (1 - \kappa)$ percent point of its bootstrap distribution. Let $B$ denote the number of bootstrap samples. The goodness-of-fit test is then carried out as follows:
\begin{enumerate}[label={\arabic*.}]
    \item For each $b \in \{1, \dots, B\}$, event times $(T_i^{1,b}, T_i^{2,b}), i = 1, \dots, n$ are generated according to model \eqref{eq:Regression Model rewritten} with $K=2$ using the estimated parameters $(\hat{\gamma}, \hat{\eta})$, where errors $(\epsilon_{1, i}, \epsilon_{2, i})$ are simulated according to a bivariate normal distribution based on $\hat{\sigma}_1$, $\hat{\sigma}_2$ and $\hat{\rho}_{12}$.
    \item We generate independent censoring times $C_i^b, i = 1, \dots, n$ by first estimating the distribution of $C$, denoted by $G_{C}$, based on a Kaplan--Meier estimator with survival times $Y$ and censoring indicators $\Delta_C$. Values for $C_i^b$ are then generated as $C_i^b = \argmax_{\tilde{t}}(\hat{G}_{C}(\tilde{t}) < U_i)$, where $U_i$ follows a standard uniform distribution.
    \item The bootstrapped times and censoring indicators are obtained by defining $Y_i^b = \min(T_i^{1,b}, T_i^{2,b}, C_i^b)$, $\Delta_{1,i}^b = I(Y_i^b = T_i^{1,b})$ and $\Delta_{2,i}^b = I(Y_i^b = T_i^{2,b})$, for all $i = 1, \dots, n$.
    \item For each bootstrap sample $(Y_i^b, \Delta_{1,i}^b, \Delta_{2,i}^b, W_i, Z_i)_{i = 1, \dots, n}$ the weighted Cramér--Von Mises type test statistic $T_{CM, b}^*$ is computed. Denote with $q^*_{CM, 1 - \kappa}$ the $(1-\kappa)$-quantile of the vector $(T^*_{CM, b})_{b = 1, \dots, B}$ of these test statistics.
    \item Reject $H_0$ on a $(1-\kappa)$-confidence level if $T_{CM} > q^*_{CM, 1 - \kappa}$.
\end{enumerate}
The performance of this test is studied in Section \ref{section: Simulation study}. In general, the test shows a good control of the type-I error and reasonable power to reject misspecified models.

\section{Simulation study}\label{section: Simulation study}
In line with Section \ref{section: Goodness of fit test}, we will consider $K = 2$ throughout the entire simulation study to keep the required computation time manageable.
\subsection{Finite sample performance}
\label{subsecsim}
Let $\tilde{X} \sim \mathcal{N}(0,1)$, $\tilde{W} \sim \text{Bernoulli}(0.5)$ and $\nu \sim \text{Logistic}(0,1)$. The endogenous variable $Z$ can then be defined as $Z=\mathbbm{1}(W^{\top}\gamma-\nu>0)$ with $\gamma=(-1,0.6,2.3)^{\top}$. Further, we construct $V$ such that:
\begin{align*}
V & = (1-Z)\bigl[ \bigl(1+\exp\{W^{\top}\gamma\}\bigl) \log \bigl(1+\exp \{W^{\top}\gamma\}\bigl)-W^{\top}\gamma \exp \{W^{\top} \gamma\} \bigl] \\
& \quad -Z \bigl[ \bigl(1+\exp \{-W^{\top}\gamma\} \bigl) \log \bigl(1+\exp \{-W^{\top}\gamma\} \bigl) + W^{\top} \gamma \exp \{-W^{\top}\gamma\}\bigl].
\end{align*}
Finally, $T^1$ and $T^2$ are constructed as:
\begin{equation*}
    \begin{array}{ll}
        T^1 = \Lambda_{\theta_1}^{-1}(X^{\top} \beta_1 + Z\alpha_1 + V\lambda_1 + \epsilon_1), \\
        T^2 = \Lambda_{\theta_2}^{-1}(X^{\top} \beta_2 + Z\alpha_2 + V\lambda_2 + \epsilon_2),
    \end{array}
\end{equation*}
for which the true parameters are $\theta_1=1$, $\theta_2=0.5$, $(\beta_{1,0},\beta_{1,1},\alpha_1,\lambda_1) =(2.5,2.6,1.8,2)$ and $(\beta_{2,0},\beta_{2,1},\alpha_2,\lambda_2)=(1.8,0.9,0.5,-2.2)$, where $\beta_{j,k}$ denotes the $k$-th element of $\beta_j$. Furthermore, it is assumed that $C \sim U[0,8]$ is generated independently from all other variables. Using these parameter values, on average $40\%$ of the generated follow-up times correspond to $T^1$, another $40\%$ corresponds to $T^2$ and the final $20\%$ to $C$. 

Three different sample sizes will be considered, namely $n=250$, $n=500$ and $n=1000$. For each sample size, $2500$ simulations are run. Interest is also in comparing the \emph{two-step estimator} with the performance of three other estimators. Firstly, we will consider an \emph{independent estimator}, which assumes that the variable $T^2$ is independent of $T^1$ given the covariates. Secondly, we consider a \emph{naive estimator} that ignores the endogeneity of $Z$ and therefore does not include $V$ in the model. The third estimator that will be considered is an \emph{oracle estimator} that assumes the control function $V$ to be observed. The results are compared in terms of bias, empirical standard deviation (ESD), root mean squared error (RMSE) and coverage rate (CR). To clarify the definition of ESD and RMSE, take $\alpha_1$ as an example. In this case, the formulas for ESD and RMSE are given by:
\begin{align*}
&\text{ESD}=\sqrt{(N-1)^{-1}\sum_{j=1}^{N}{\bigl[ (\hat{\alpha}_1)_j-\bar{\alpha}_1\bigl]^2}}, \quad \text{with } \bar{\alpha}_1 = N^{-1}\sum_{k=1}^{N}{(\hat{\alpha}_1)_k}, \\
&\text{RMSE}=\sqrt{N^{-1}\sum_{j=1}^{N}{\bigl[ (\hat{\alpha}_1)_j-\alpha_1^*\bigl]^2}}, \quad \text{with } \alpha_1^* \text{ the true parameter value.}
\end{align*}
Here $(\hat{\alpha}_1)_j$ denotes the estimate of $\alpha_1^*$ in the $j$-th simulation and $N$ denotes the number of simulations. When the bias decreases, it holds that $\bar{\alpha}_1$ converges to $\alpha_1^*$ and hence the ESD and RMSE converge to the same value. The coverage rate is calculated as the percentage of simulations in which the true parameter is contained within the $95\%$ confidence interval with the variance computed as in Section \ref{section: Asymptotic variance}. The results of the simulations can be found in Table \ref{tab: results simulation design 4}. Three other simulation settings were considered as well, which show similar results and can be found in Section D of the Supplementary Material. The results show that the \emph{naive estimator} has a large bias for most of the parameters and that this bias does not shrink towards zero as the sample size increases. Besides, the coverage rate is very small for most parameters, meaning that the confidence intervals rarely contain the true parameter value. Hence, ignoring the endogeneity of a covariate leads to poor results. Secondly, we can observe similar but less extreme results using the \emph{independent estimator}. Although the bias does not decrease towards zero either, it is seen to be lower than using the \emph{naive estimator}. For the \emph{two-step estimator}, it can be seen that the bias decreases towards zero as the sample size increases. Besides, when the sample size is $1000$, the coverage rate is close to $95\%$ for most parameters, indicating that the estimated standard errors are asymptotically valid. Although the coverage rate is still considerably lower than $95\%$ for some parameters (e.g. $\alpha_2$), additional simulations with an even larger sample size showed that these coverage rates converge to $95\%$ as well. Next, we can see that the ESD and the RMSE converge to the same decreasing value as the sample size increases. Lastly, comparing the results of bias and coverage of the \emph{two-step estimator} to the results of the \emph{oracle estimator}, it can be concluded that they are fairly similar. This implies that the error made by estimating the control function $V$ is much smaller than the error resulting from the second step estimation.

\begin{table}[!ht]
\centering
\scalebox{0.88}{
\begin{tabular}{|c|rrrr|rrrr|rrrr|}
  \hline \multicolumn{5}{|c}{$n = 250$} & \multicolumn{4}{|c}{$n = 500$}& \multicolumn{4}{|c|}{$n = 1000$} \\ \hline \multicolumn{13}{|c|}{two-step estimator} \\  \hline
 & Bias & ESD & RMSE & CR & Bias & ESD & RMSE & CR & Bias & ESD & RMSE & CR \\ 
  \hline
$\beta_{1,0}$ & -0.040 & 0.594 & 0.595 & 0.944 & -0.031 & 0.405 & 0.406 & 0.954 & -0.020 & 0.286 & 0.287 & 0.950 \\ 
  $\beta_{1,1}$ & -0.012 & 0.263 & 0.263 & 0.943 & -0.004 & 0.180 & 0.180 & 0.948 & -0.000 & 0.127 & 0.127 & 0.951 \\ 
  $\alpha_1$ & 0.029 & 0.839 & 0.839 & 0.948 & 0.026 & 0.588 & 0.589 & 0.949 & 0.019 & 0.413 & 0.413 & 0.954 \\ 
  $\lambda_1$ & -0.006 & 0.330 & 0.330 & 0.929 & 0.005 & 0.231 & 0.231 & 0.942 & 0.005 & 0.163 & 0.163 & 0.944 \\ 
  $\beta_{2,0}$ & 0.003 & 0.502 & 0.501 & 0.932 & 0.006 & 0.356 & 0.356 & 0.935 & 0.005 & 0.249 & 0.249 & 0.924 \\ 
  $\beta_{2,1}$ & 0.007 & 0.269 & 0.269 & 0.948 & 0.006 & 0.190 & 0.190 & 0.949 & -0.002 & 0.137 & 0.137 & 0.940 \\ 
  $\alpha_2$ & 0.035 & 1.098 & 1.098 & 0.914 & -0.002 & 0.674 & 0.673 & 0.917 & -0.001 & 0.476 & 0.476 & 0.907 \\ 
  $\lambda_2$ & -0.006 & 0.391 & 0.391 & 0.942 & -0.004 & 0.275 & 0.275 & 0.936 & -0.003 & 0.192 & 0.192 & 0.922 \\ 
  $\sigma_1$ & -0.023 & 0.080 & 0.083 & 0.914 & -0.012 & 0.056 & 0.057 & 0.935 & -0.006 & 0.039 & 0.039 & 0.942 \\ 
  $\sigma_2$ & -0.022 & 0.104 & 0.106 & 0.915 & -0.011 & 0.071 & 0.072 & 0.936 & -0.005 & 0.051 & 0.051 & 0.943 \\ 
  $\rho_{12}$ & 0.014 & 0.164 & 0.165 & 0.879 & 0.007 & 0.103 & 0.103 & 0.914 & 0.003 & 0.070 & 0.070 & 0.939 \\ 
  $\theta_1$ & -0.006 & 0.045 & 0.045 & 0.925 & -0.003 & 0.031 & 0.031 & 0.940 & -0.002 & 0.022 & 0.022 & 0.934 \\ 
  $\theta_2$ & -0.001 & 0.085 & 0.085 & 0.929 & -0.001 & 0.058 & 0.058 & 0.946 & -0.001 & 0.041 & 0.041 & 0.945 \\ 
   \hline \multicolumn{13}{|c|}{naive estimator} \\ \hline $\beta_{1,0}$ & 2.752 & 1.013 & 2.933 & 0.062 & 2.733 & 0.697 & 2.822 & 0.008 & 2.721 & 0.522 & 2.771 & 0.005 \\ 
  $\beta_{1,1}$ & 0.607 & 0.282 & 0.669 & 0.279 & 0.610 & 0.203 & 0.643 & 0.048 & 0.613 & 0.158 & 0.633 & 0.002 \\ 
  $\alpha_1$ & -4.730 & 0.921 & 4.819 & 0.003 & -4.716 & 0.607 & 4.756 & 0.002 & -4.702 & 0.431 & 4.722 & 0.000 \\ 
  $\beta_{2,0}$ & -2.367 & 0.216 & 2.377 & 0.000 & -2.369 & 0.176 & 2.376 & 0.000 & -2.363 & 0.152 & 2.368 & 0.000 \\ 
  $\beta_{2,1}$ & -0.583 & 0.208 & 0.619 & 0.160 & -0.581 & 0.170 & 0.606 & 0.010 & -0.587 & 0.145 & 0.604 & 0.000 \\ 
  $\alpha_2$ & 5.517 & 1.022 & 5.611 & 0.008 & 5.452 & 0.666 & 5.494 & 0.004 & 5.428 & 0.579 & 5.459 & 0.002 \\ 
  $\sigma_1$ & 0.609 & 0.142 & 0.626 & 0.001 & 0.622 & 0.100 & 0.630 & 0.000 & 0.629 & 0.076 & 0.634 & 0.000 \\ 
  $\sigma_2$ & 0.511 & 0.158 & 0.535 & 0.024 & 0.516 & 0.110 & 0.528 & 0.001 & 0.520 & 0.083 & 0.526 & 0.000 \\ 
  $\rho_{12}$ & -0.526 & 0.434 & 0.681 & 0.642 & -0.495 & 0.274 & 0.566 & 0.246 & -0.483 & 0.189 & 0.519 & 0.029 \\ 
  $\theta_1$ & 0.026 & 0.063 & 0.068 & 0.904 & 0.032 & 0.045 & 0.056 & 0.845 & 0.035 & 0.033 & 0.048 & 0.735 \\ 
  $\theta_2$ & 0.194 & 0.086 & 0.212 & 0.396 & 0.194 & 0.058 & 0.202 & 0.117 & 0.195 & 0.041 & 0.199 & 0.005 \\ 
   \hline \multicolumn{13}{|c|}{independent estimator} \\ \hline $\beta_{1,0}$ & 0.696 & 0.578 & 0.904 & 0.711 & 0.692 & 0.399 & 0.799 & 0.551 & 0.689 & 0.283 & 0.745 & 0.331 \\ 
  $\beta_{1,1}$ & 0.119 & 0.249 & 0.276 & 0.902 & 0.121 & 0.171 & 0.209 & 0.877 & 0.123 & 0.120 & 0.171 & 0.808 \\ 
  $\alpha_1$ & -0.552 & 0.858 & 1.021 & 0.838 & -0.552 & 0.599 & 0.815 & 0.789 & -0.547 & 0.423 & 0.691 & 0.691 \\ 
  $\lambda_1$ & 0.061 & 0.337 & 0.342 & 0.952 & 0.068 & 0.235 & 0.244 & 0.958 & 0.069 & 0.166 & 0.179 & 0.948 \\ 
  $\beta_{2,0}$ & 0.096 & 0.517 & 0.526 & 0.949 & 0.097 & 0.366 & 0.379 & 0.942 & 0.098 & 0.259 & 0.276 & 0.926 \\ 
  $\beta_{2,1}$ & -0.087 & 0.264 & 0.278 & 0.920 & -0.085 & 0.185 & 0.204 & 0.910 & -0.093 & 0.133 & 0.162 & 0.874 \\ 
  $\alpha_2$ & 0.299 & 1.097 & 1.137 & 0.902 & 0.266 & 0.702 & 0.751 & 0.892 & 0.256 & 0.501 & 0.563 & 0.871 \\ 
  $\lambda_2$ & -0.078 & 0.403 & 0.410 & 0.953 & -0.074 & 0.283 & 0.292 & 0.940 & -0.073 & 0.199 & 0.212 & 0.924 \\ 
  $\sigma_1$ & -0.019 & 0.084 & 0.086 & 0.925 & -0.005 & 0.059 & 0.059 & 0.940 & 0.001 & 0.041 & 0.041 & 0.949 \\ 
  $\sigma_2$ & 0.009 & 0.109 & 0.110 & 0.933 & 0.020 & 0.076 & 0.078 & 0.935 & 0.026 & 0.054 & 0.060 & 0.920 \\ 
  $\theta_1$ & -0.005 & 0.045 & 0.045 & 0.931 & -0.004 & 0.031 & 0.031 & 0.942 & -0.004 & 0.022 & 0.023 & 0.932 \\ 
  $\theta_2$ & -0.022 & 0.087 & 0.090 & 0.926 & -0.021 & 0.060 & 0.064 & 0.931 & -0.021 & 0.042 & 0.047 & 0.918 \\ 
   \hline \multicolumn{13}{|c|}{oracle estimator} \\ \hline $\beta_{1,0}$ & -0.007 & 0.409 & 0.409 & 0.926 & -0.009 & 0.270 & 0.270 & 0.945 & -0.009 & 0.188 & 0.188 & 0.946 \\ 
  $\beta_{1,1}$ & -0.005 & 0.152 & 0.152 & 0.938 & -0.002 & 0.103 & 0.103 & 0.945 & -0.002 & 0.071 & 0.071 & 0.948 \\ 
  $\alpha_1$ & -0.010 & 0.485 & 0.485 & 0.929 & -0.006 & 0.326 & 0.326 & 0.949 & 0.001 & 0.227 & 0.227 & 0.952 \\ 
  $\lambda_1$ & -0.008 & 0.178 & 0.178 & 0.934 & -0.004 & 0.120 & 0.120 & 0.948 & -0.002 & 0.084 & 0.084 & 0.950 \\ 
  $\beta_{2,0}$ & -0.000 & 0.290 & 0.290 & 0.930 & -0.005 & 0.198 & 0.198 & 0.940 & -0.001 & 0.143 & 0.143 & 0.928 \\ 
  $\beta_{2,1}$ & 0.003 & 0.154 & 0.154 & 0.940 & 0.007 & 0.106 & 0.106 & 0.945 & 0.001 & 0.077 & 0.077 & 0.939 \\ 
  $\alpha_2$ & 0.070 & 1.007 & 1.009 & 0.921 & 0.016 & 0.415 & 0.415 & 0.921 & 0.008 & 0.300 & 0.300 & 0.918 \\ 
  $\lambda_2$ & -0.005 & 0.245 & 0.245 & 0.932 & 0.001 & 0.166 & 0.166 & 0.933 & 0.000 & 0.120 & 0.120 & 0.932 \\ 
  $\sigma_1$ & -0.024 & 0.080 & 0.083 & 0.927 & -0.012 & 0.056 & 0.057 & 0.944 & -0.006 & 0.039 & 0.039 & 0.943 \\ 
  $\sigma_2$ & -0.022 & 0.104 & 0.106 & 0.928 & -0.011 & 0.071 & 0.072 & 0.939 & -0.005 & 0.051 & 0.051 & 0.943 \\ 
  $\rho_{12}$ & 0.013 & 0.157 & 0.158 & 0.926 & 0.009 & 0.095 & 0.095 & 0.933 & 0.004 & 0.062 & 0.063 & 0.948 \\ 
  $\theta_1$ & -0.005 & 0.043 & 0.043 & 0.942 & -0.003 & 0.030 & 0.030 & 0.946 & -0.002 & 0.021 & 0.021 & 0.939 \\ 
  $\theta_2$ & -0.004 & 0.084 & 0.085 & 0.941 & -0.002 & 0.058 & 0.058 & 0.952 & -0.002 & 0.041 & 0.041 & 0.942 \\ 
   \hline
\end{tabular}
}
\caption{Estimation results for 2500 simulations where 40\% of the observations correspond to $T_1$, another $40\%$ correspond to $T_2$ and the final $20\%$ correspond to $C$. Given are the bias, empirical standard deviation (ESD), root mean squared error (RMSE) and coverage rate (CR).}
 \label{tab: results simulation design 4}
\end{table}

\subsection{Model under misspecification} \label{section: model under misspecification}
We now study the performance of the model when either the control function is misspecified or the error terms are not normally distributed, even after applying the Yeo--Johnson transformation. In the following, we always used $3$ different sample sizes (namely $250, 500$ and $1000$) and $500$ simulated data sets for each sample size. The results of these simulations can be found in Section E.2 of the Supplementary Material.

\subsubsection{Misspecified control function}
Consider the same setting as in Section \ref{subsecsim}, with a binary variable $Z$. To identify the control function $V$, we model $Z = \mathbbm{1}(W^\top \gamma - \nu > 0)$, for some specified distribution of $\nu$. In Section \ref{subsecsim}, we assumed $\nu \sim \text{Logistic}(0, 1)$. However, other distributions could have been assumed on $\nu$ as well. In particular, two common choices for the link function are the probit link and the cumulative log-log link. To investigate the effect of possible misspecification, data will be simulated for which $\nu$ follows a standard normal distribution (for the probit link) or $\nu$ follows a standard Gumbel distribution (cumulative log-log link). The assumed model will always make use of the logit link and hence the control function is misspecified. The necessary formulas and derivations can be found in Section E.1 of the Supplementary Material. 

We find that for both considered distributions on $\nu$, the estimates of $\lambda_1$ and $\lambda_2$ are biased and the theoretical confidence intervals almost never contain the true parameter values. However, if $\nu$ follows a normal distribution (probit link), the bias of the estimated causal effect of $Z$ on $T^1$ seems to be relatively low and the coverage rate is still reasonably high, even though $\lambda_1 V$ does not completely capture the part of $u_1$ that is correlated with $Z$. Also for the other parameters, the bias and coverage rate are still acceptable in our opinion. When $\nu$ follows a standard Gumbel distribution, the results become less precise. The bias of $\alpha_1$ has increased considerably. This could be expected as the logit link is less similar to the cumulative log-log link than to the probit link.

\subsubsection{Errors deviating from normal distribution} \label{section: Errors deviating from normal distribution}
Next, it is investigated what happens when the error terms are not normally distributed. To this end, the performance of the model will be evaluated when the error terms follow a bivariate skew-normal distribution or a bivariate t-distribution (with $3$ degrees of freedom). In this way, the effect of asymmetry or higher-than-assumed kurtosis of the error distribution can be studied. To simulate the bivariate skew-normal distribution, the ideas of \cite{Azzalini1996} are used. The index of skewness of the marginal distribution of the error terms is chosen to be $0.92$. Besides, the performance of the model when the error terms are heteroscedastic will also be investigated. We consider all these misspecifications in the case where both $\tilde{W}$ and $Z$ are continuous. To generate the data, assume that $\tilde{X} \sim \mathcal{N}(0,1), \tilde{W} \sim U[0,2]$,  $\nu \sim \mathcal{N}(0,2)$ and $Z=W^{\top}\gamma+\nu$. The parameter values are still the same as in Section \ref{subsecsim}.

In summary, we found that for a skew-normal or t-distribution, the impact on the bias as well as coverage rate of $\alpha_1$, which is the parameter of main interest, is rather small. In the case of heteroscedastic errors, we found that the bias increases with increasing heteroskedasticity. However, for all considered cases of misspecified error distributions, the impact on the estimation of $\alpha_1$ seems to be relatively limited.

\subsection{Goodness-of-fit test}
We also study the type-I error and power of the goodness-of-fit test as described in Section \ref{section: Goodness of fit test}. To investigate the type-I error, we can repeatedly simulate data sets according to model \eqref{eq:Regression Model rewritten} with specified coefficients and refit the model to the generated data. We then count the number of times the fitted model is rejected and from this, we estimate the type-I error since we know that every rejection is a false rejection. Similarly, we can gauge the power of the test by repeatedly simulating data sets according to the misspecified models discussed in Section \ref{section: model under misspecification}.

We perform this simulation study using $500$ simulated data sets of size $n = 1000$ and $n = 2000$. The distribution of the test statistic will be based on $B = 250$ and $B = 500$ bootstrap samples when assessing the type-I error of the test and $B = 250$ bootstrap samples when assessing its power. In the following, we will refer to the scenario where the data are generated according to model \eqref{eq:Regression Model rewritten} as scenario $0$. The cases where the control function is based on a probit link or cumulative log-log link are referred to as scenarios $1$-a and $1$-b, respectively. The case in which the data is generated using the skew-normal distribution of \cite{Azzalini1996} and a bivariate t-distribution will be referred to as scenario $2$-a and $2$-b, respectively. Lastly, the case wherein the errors are heteroscedastic will be called scenario $2$-c. We will use the same parameters as those chosen in Section \ref{section: model under misspecification} for all of these misspecified distributions. Furthermore, we assume that $Z$ and $W$ are continuous in all scenarios except in scenarios $1$ -a and $1$ -b, where we have to let $Z$ be a binary variable in order to investigate the effect of misspecification of the link function.

\begin{table}[t]
    \centering
    \begin{tabular}{cccccc}
        \hline
        \hline
        & \multirow{2}{*}{$B$} & \multirow{2}{*}{Sample size} & \multirow{2}{*}{Scenario} & \multicolumn{2}{c}{Rejection rate at}\\
        \cmidrule(rl){5-6}
        & & & & $5\%$ & $10\%$\\
        \cmidrule(rl){1-6}
        \multirow{4}{*}{Type-I error} & $250$ & $n = 1000$ & 0 & 0.042 & 0.082\\
        & $250$ & $n = 2000$ & 0 & 0.044 & 0.098\\
        & $500$ & $n = 1000$ & 0 & 0.044 & 0.084\\
        & $500$ & $n = 2000$ & 0 & 0.044 & 0.096\\
        \cmidrule(rl){1-6}
        \multirow{10}{*}{Power} & 250 & $n = 1000$ & 1-a & 0.046 & 0.094\\
        & 250 & $n = 1000$ & 1-b & 0.278 & 0.426\\
        & 250 & $n = 1000$ & 2-a & 0.050 & 0.100\\
        & 250 & $n = 1000$ & 2-b & 0.172 & 0.220\\
        & 250 & $n = 1000$ & 2-c & 0.198 & 0.304\\

        \cmidrule(rl){2-6}
        & 250 & $n = 2000$ & 1-a & 0.030 & 0.076\\
        & 250 & $n = 2000$ & 1-b & 0.578 & 0.726\\
        & 250 & $n = 2000$ & 2-a & 0.054 & 0.110\\
        & 250 & $n = 2000$ & 2-b & 0.344 & 0.424\\
        & 250 & $n = 2000$ & 2-c & 0.284 & 0.390\\
        \hline
    \end{tabular}
    \caption{Simulation results for 500 simulations of the power and type-I error of the goodness-of-fit test. Scenario 0 corresponds to correctly specifying the model, misspecification of the control function based on a probit or cumulative log-log link are scenarios 1-a and 1-b respectively. Scenarios 2-a and 2-b correspond to generating the data from a skew-normal and bivariate t-distribution respectively. Lastly, scenario 2-c deals with heteroscedastic errors.}
    \label{tab: rejection rates}
\end{table}

The results are shown in Table \ref{tab: rejection rates}. From this table, it can be seen that the type-I error of the test when $\kappa = 0.05$ is close to $5\%$ albeit slightly conservative. When the errors come from a bivariate t-distribution or are heteroscedastic (i.e. scenarios $2$-b and $2$-c), we can see that the test has some power to detect that misspecification and that this power increases with the sample size. It could be remarked that the rejection rates in these cases are not very large. This does not necessarily imply that our proposed test lacks power, but it could be the result of the extra flexibility introduced by applying the Yeo--Johnson transformation to the observed times. After all, it has the effect of making a random variable more normally distributed. We can also see the results of this effect when the errors are distributed according to a skew-normal distribution (scenario $2$-a), where the applied transformation seemingly completely counteracts the misspecification of the error distribution. Furthermore, we can observe something similar for the power of the test to detect a misspecification of the control function. When the control function is incorrectly based on a logit link instead of the probit link (scenario $1$-a), the test is unable to detect it. This seems understandable since both link functions are similar. The strange decrease in power when the sample size increases is likely due to noise. On the contrary, the test has sufficient power to detect that the control function is misspecified when it should have been based on the cumulative log-log link (scenario $1$-b). Moreover, we can see that the power of the test to detect such a misspecification grows substantially with increasing sample size.

At first glance, one could conclude from this simulation that the goodness-of-fit test lacks power. However, some considerations should be taken into account. Firstly, it is important to note the difficulty of the simulation settings, which often pertain to only a slight misspecification in an otherwise fully correctly specified model. Secondly, we reiterate that the test is conservative, since we can only test the fit of the model based on the observed data distribution (cf. Section \ref{section: Goodness of fit test}). Therefore, we deem these results to be satisfactory.


\section{Data application}\label{section:dataapp}
The modeling approach as discussed in this paper will now be applied to estimate the effect of Job Training Partnership Act (JTPA) services on time until employment. The data come from a large-scale study that was commissioned by the US Department of Labor in 1986, meant to assess the benefits and cost of several training programs related to job employability \citep{Bloom1997}.

The study was designed to estimate the effects of the training programs for several key demographic groups. To this end, participants were assigned to either a control or a treatment group. Because the staff at the Service Delivery Areas did not always adhere to the randomization rules, some two-sided noncompliance took place. More precisely, $3\%$ of the control group members still participated in the JTPA services. The focus will be on the $212$ married, white men without children who did not have a job at the time they were randomized. We use this relatively small stratum to limit the computational intensity of the goodness-of-fit test later on. The exogenous predictors in our model will be the participant's age and a variable indicating whether the participant has achieved a high school diploma or GED (\textit{hsged}). The endogenous variable $Z$ is an indicator of whether the participant followed a training program and its instrument $\Tilde{W}$ is a binary variable indicating if the participant was assigned to the control or treatment group ($0$ and $1$ respectively). Because of the two-sided noncompliance, $Z$ is confounded as individuals moved themselves between the control and treatment arm in a non-random way. We argue that $\Tilde{W}$ is an appropriate instrument for $Z$ as it is randomly assigned (conditional on the measured covariates) and hence uncorrelated with the error terms $(\epsilon_1,\epsilon_2)$, clearly correlated with $Z$ and it can only affect the time until employment through $Z$. We will assume that $\nu$ follows a standard logistic distribution. 

Our interest lies in the time between randomization and employment $(T_1)$. To measure this time, researchers invited the participants to one or possibly two follow-up interviews. For the participants who only received an invitation for the first interview, $T_1$ is recorded precisely if they are employed at the time of the interview. In the other case, the observation is regarded as lost to follow-up at the time of the interview $(T_2)$. For individuals who were invited to the second interview and participated in it, the time until employment is recorded precisely if the individual is employed at the time of this second interview. Otherwise, the observation is regarded as lost to follow-up at the time of this interview. The individuals who were invited but did not participate in the second survey are considered lost to follow-up at the time of the first interview unless they were already employed at the time of the first survey. In the latter case, their time until employment was already recorded precisely. By defining $T_1$ and $T_2$ in this way, it is likely that $T_1$ and $T_2$ are dependent. This is because the decision to participate in the second interview could be influenced by whether or not the participant has found a job between the first and second interview. In this way, $13\%$ of the observations are censored.

\begin{figure}[H]
    \centering
    \includegraphics[width=0.8\linewidth]{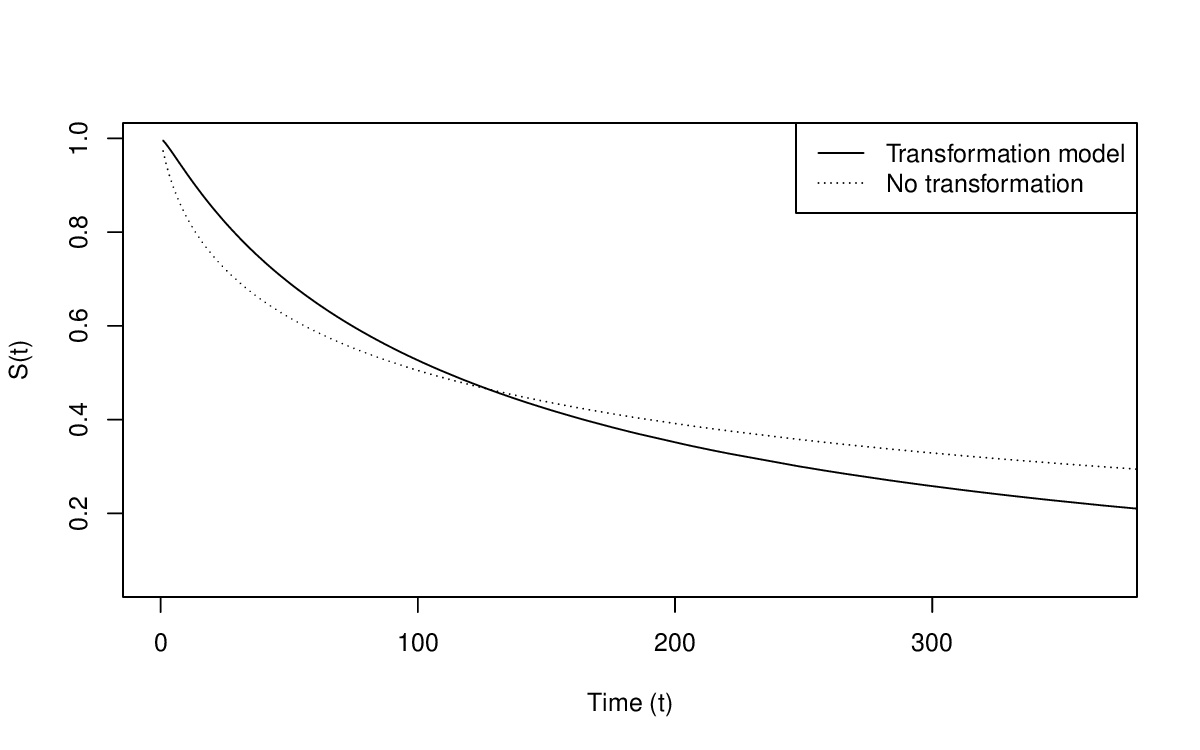}
    \includegraphics[width=0.8\linewidth]{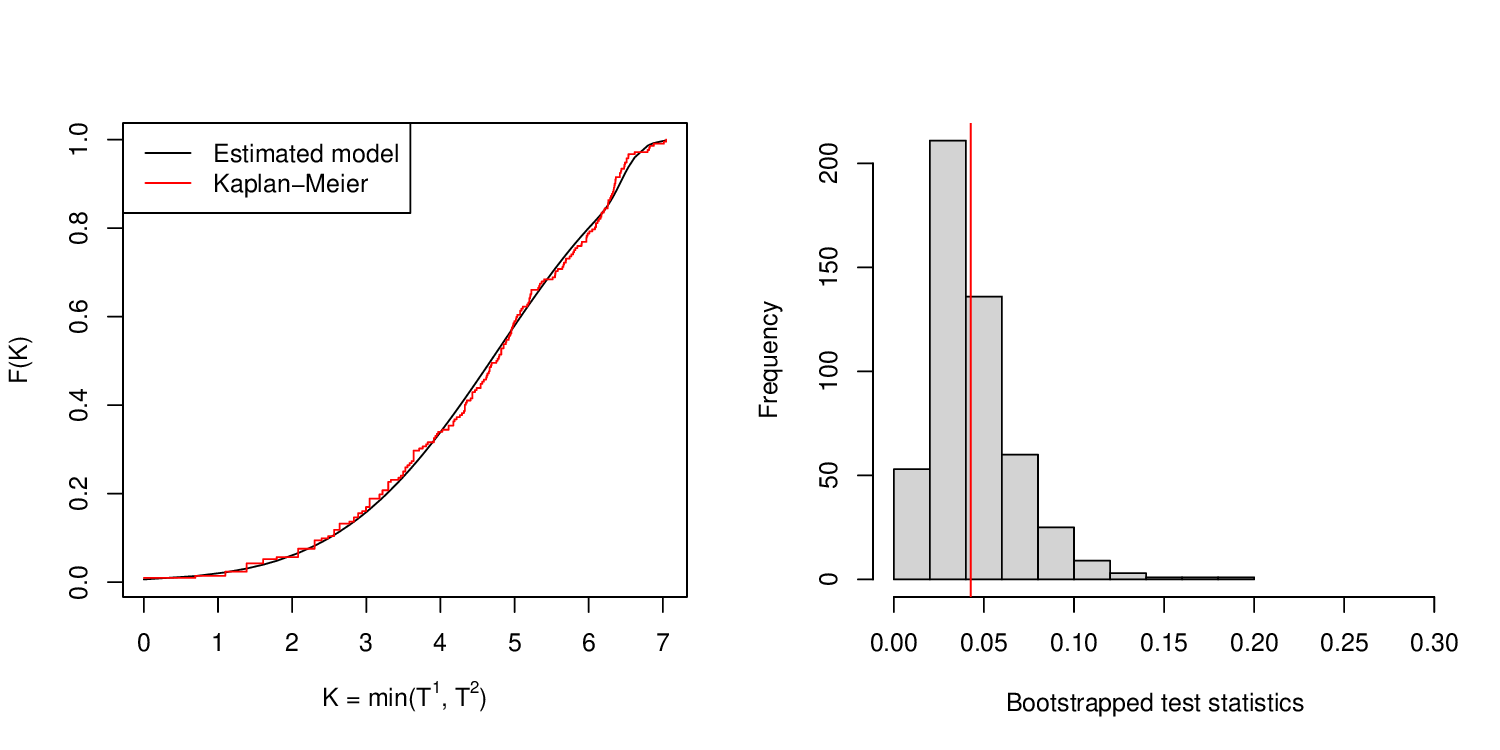}
    \caption{Top panel: the estimated survival curves using the two-step estimator with (solid curve) and without (dotted curve) transforming the response. The bottom panels relate to the goodness of fit test. The bottom left panel compares the distributions of $K = \min(T^1, T^2)$ implied by the model and based on the Kaplan-Meier estimator. The bottom right panel shows the bootstrap distribution of the goodness-of-fit test statistic. The vertical red line represents the observed value of the test statistic.}
    \label{fig:Comparison_CVK_YJ}
\end{figure}

In this empirical application, three estimators are compared: the \emph{two-step estimator}, the \textit{naive estimator} and the \textit{independent estimator}. The estimated coefficients, alongside their standard errors and $p$-values, are shown in Table \ref{tab:DataApplication}. Note that for the transformation parameters $\theta_1$ and $\theta_2$, the $p$-value measures the significance of the difference with respect to unity. An application of the goodness-of-fit test for our proposed model based on $500$ bootstrap samples gives a $p$-value of $0.424$, hence we do not reject that the proposed model fits the given data set. The bootstrap distribution of the goodness-of-fit test statistic is plotted in Figure \ref{fig:Comparison_CVK_YJ}.

From Table \ref{tab:DataApplication} it can be seen that all three models use Yeo-Johnson transformations with parameters significantly different from unity in order to improve normality. The main conclusion from this study is that the proposed model estimates the effect of the treatment to be positive. That is, participants who attend the job training programs generally take longer to find a job than participants who do not attend the programs. We believe that this effect is due to the programs having little to no influence on the employability of the participants, and participants who take longer to find a job are self-selecting into the job training programs. We do remark that this conclusion is specific to the stratum under investigation. For example, when considering married, black fathers, no such effect is detected.

Lastly, we plot the estimated survival curves using the \textit{two-step estimator} with (solid curve) and without (dotted curve) applying a Yeo--Johnson transformation. We do this for a man who is $30$ years old, has not obtained a high school diploma, was assigned to the treatment group and complied. It can be seen that the curves differ substantially. However, the estimated median survival times are close: the transformation model estimates it at $111$ days, while the model assuming $\theta_1= \theta_2 = 1$ estimates it to be at $103$ days.

\begin{table}[!ht]
    \centering
    \begin{adjustbox}{width=0.95\textwidth}
    \begin{tabular}{|r|rrr|rrr|rrr|}
        \hline
        & \multicolumn{3}{c|}{two-step estimator} & \multicolumn{3}{c|}{naive estimator} & \multicolumn{3}{c|}{independent estimator}\\
        \hline
          & Estimate & SE & $p$-value & Estimate & SE & $p$-value & Estimate & SE & $p$-value \\ 
        \hline
        $\beta_{1,0}$            &  4.85 & 0.17 & 0.00 &  5.17 & 0.16 & 0.00 &  3.98 & 0.86 & 0.00  \\ 
        $\beta_{1,\text{age}}$   &  0.03 & 0.01 & 0.04 &  0.03 & 0.01 & 0.01 &  0.03 & 0.01 & 0.01   \\ 
        $\beta_{1,\text{hsged}}$ & -0.41 & 0.24 & 0.09 & -0.38 & 0.34 & 0.25 & -0.28 & 0.31 & 0.36   \\ 
        $\alpha_1$               &  0.69 & 0.19 & 0.00 &  0.21 & 0.35 & 0.54 &  0.74 & 0.70 & 0.29  \\ 
        $\lambda_1$              &  0.19 & 0.07 & 0.01 &       &      &      &  0.23 & 0.26 & 0.38  \\
        \hline

        $\beta_{2,0}$            &  2.87 & 0.26 & 0.00 &  3.50 & 0.01 & 0.00 &  7.66 & 0.02 & 0.00  \\ 
        $\beta_{2,\text{age}}$   &  0.00 & 0.00 & 0.13 &  0.00 & 0.00 & 0.40 &  0.00 & 0.00 & 0.89 \\ 
        $\beta_{2,\text{hsged}}$ & -0.02 & 0.00 & 0.00 & -0.03 & 0.03 & 0.34 & -0.08 & 0.09 & 0.38  \\ 
        $\alpha_2$               & -0.01 & 0.01 & 0.26 &  0.02 & 0.03 & 0.45 & -0.36 & 0.24 & 0.12  \\ 
        $\lambda_2$              & -0.01 & 0.00 & 0.15 &       &      &      & -0.17 & 0.09 & 0.07 \\ 
        \hline
        $\sigma_1$ & 2.41 & 0.31 & 0.00 & 2.45 & 0.24 & 0.00 & 2.02 & 0.51 & 0.00 \\ 
        $\sigma_2$ & 0.11 & 0.02 & 0.00 & 0.17 & 0.03 & 0.00 & 0.25 & 0.04 & 0.00 \\ 
        $\rho_{12}$     & 0.98 & 0.02 & 0.00 & 0.98 & 0.02 & 0.00 &      &      &         \\ 
        $\theta_1$ & 1.25 & 0.07 & 0.00 & 1.26 & 0.06 & 0.00 & 1.12 & 0.16 & 0.45 \\ 
        $\theta_2$ & 0.38 & 0.08 & 0.00 & 0.57 & 0.02 & 0.00 & 1.10 & 0.01 & 0.00 \\ 
        \hline
    \end{tabular}
    \end{adjustbox}
    \caption{Estimation results for the two-step, naive and independent estimator. Given are the parameter estimate, estimated standard error (SE) and $p$-value.}
    \label{tab:DataApplication}
\end{table}

\subsection*{Funding}
I. Van Keilegom gratefully acknowledges funding from the FWO and F.R.S.-FNRS under the Excellence of Science (EOS) programme, project ASTeRISK (grant No. 40007517), and from the FWO (senior research projects fundamental research, grant no. G047524N). G. Crommen is funded by a PhD fellowship from the Research Foundation - Flanders (grant number 11PKA24N).


\bibliography{references}

\end{document}